\magnification 1200

\font\Bbb=msbm10
\def\BBB#1{\hbox{\Bbb#1}}

\font\Frak=eufm10
\def\frak#1{{\hbox{\Frak#1}}}

\font\small=cmr8

\def\L{{\cal L}}
\def\J{{\cal J}}
\def\F{{\cal F}}

\def\C{\BBB C}
\def\Z{\BBB Z}

\def\R{\BBB R}
\def\Q{\BBB Q}

\def\T{{\BBB T}^n}
\def\Tone{{\BBB T}^1}

\def\g{{\frak g}}
\def\wg{\widetilde \g}
\def\dg{\dot \g}

\def\Diff{{\rm Diff} (M)}
\def\Vect{{\rm Vect} (M)}
\def\ad{{\rm ad}}
\def\dim{{\rm dim \ }}
\def\End{\hbox{\rm End }}
\def\Der{\hbox{\rm Der }}
\def\Id{\hbox{\rm Id }}
\def\Span{\hbox{\rm Span }}
\def\VT{{\rm Vect} (\T)}
\def\FT{\F(\T)}

\def\u{{\overline u}}
\def\d{\partial}
\def\gl{gl_n}
\def\e#1{e^{2\pi i #1 x}}

\centerline{
\bf
Jet Modules.}

\

\centerline{
{\bf Yuly Billig}  
%\footnote{*}{Research supported by the  Natural Sciences and
%Engineering Research Council of Canada.}
\footnote{}{AMS subject classification: 17B66, 58A20.}
}

\

\

{\narrower\noindent\small{\bf Abstract}:
In this paper we classify indecomposable modules for the Lie algebra
of vector fields on a torus that admit a compatible action of the algebra
of functions. An important family of such modules is given by spaces of jets
of tensor fields.
\par}

\

\

{\bf 0. Introduction.}

 In recent years there was a substantial progress in representation 
theory of infinite-dimensional Lie algebras of rank $n>1$, toroidal
Lie algebras in particular. In this paper we turn our attention to 
another Lie algebra of rank $n$, the Lie algebra $W_n$
of vector fields on an $n$-dimensional torus $\T$.

 An important class
of irreducible representations for $W_n$ has its origin in differential
geometry -- these are the modules of tensor fields on a torus. 
In addition to being modules for the Lie algebra of vector fields,
tensor fields also admit multiplication by functions. For the torus, which is
a flat manifold, the spaces of tensor fields are free modules of a finite
rank over the commutative algebra of functions $\FT$.
 We formalize this property in the definition of a category $\J$ 
of $W_n$-modules (cf., [R2]).

 We also discuss another class of $W_n$-modules of a geometric nature
-- the modules of jets of tensor fields [S]. Jets of functions are used as
a tool for the symmetry analysis for partial differential equations [O].
The action on a space of jets of the Lie algebra of vector fields, known
under the term ``prolongation of vector fields'', plays a key role in that 
theory. 

From the algebraic point of view, 
jet modules are typically not irreducible, but are often indecomposable.
 The goal of the present paper is to classify indecomposable modules in 
category $\J$.

 Let us state our result in case $n=1$, for the sake of simplicity of 
notations.

 {\bf Theorem.} There is a 1-1 correspondence between indecomposable
$W_1$-modules $J$
in category $\J$ and pairs $(\lambda, U)$, where $\lambda \in \C/\Z$ and 
$U$ is a finite-dimensional indecomposable module for a Lie algebra
$$ W_1^+ = \Span \left< z^m {d \over dz} \; \bigg| \; m\geq 1 \right> .$$
Such a correspondence is given by the tensor product decomposition
$$ J = \F (\Tone) \otimes U ,$$
where $W_1$ acts according to the formula
$$ \left( {1\over 2\pi i} \e{s} {d\over dx}\right) \left( \e{m} \otimes u \right)
= (m+\lambda) \e{(s+m)} \otimes u 
+ \sum\limits_{b\geq 1} {s^b \over b!} \e{(s+m)} \otimes 
\rho \left( z^b {d\over dz} \right) u .$$
The sum in the right hand side is finite since for every finite-dimensional
representation of $W_n^+$, \ 
$\rho \left( z^b {d\over dz} \right) = 0$ for $b \gg 1$.

 This result reduces the classification of modules in $\J$ to a problem
in a completely finite-dimensional set-up -- describing finite-dimensional
representations of certain finite-dimensional Lie algebras.

 The case of rank one is of interest by itself since in this case we deal 
with the representations of the Virasoro algebra with a trivial action of 
the center, and the Virasoro algebra plays a prominent role in 
applications to physics. Only in the case $n=1$ we have a complete 
classification of irreducible modules with finite-dimensional weight
spaces [M]. Very little has been known about indecomposable $W_n$-modules
even in rank 1 case.

 The concept of a polynomial module, introduced in [BB] (see also [BZ]),
has turned out to be extremely useful for the present paper. We prove that
all modules in category $\J$ are polynomial modules, and this property
allows us to establish the classification result. It is interesting that
unlike all previously known examples, the degrees of the structure
polynomials for the modules in $\J$ can be arbitrarily high.

 The technique developed in this paper also allows us to recover 
Eswara Rao's classification [R2] of irreducible modules in category $\J$,
significantly simplifying his proof.

 The structure of the paper is the following. In Section 1 we discuss 
the module structure on the space of jets of tensor fields.
Motivated by the construction of jet modules, we introduce a category
$\J$ in Section 2, and state our main theorems at the end of the section.
 In Sections 3  and  4 we prove 
the classification results for the modules in category $\J$.
In the final Section 5, we give the analogous result for the semidirect
product of $W_n$ with a multi-loop algebra.

\

\

{\bf 1. Jets of tensor fields.}

\

 The group of diffeomorphisms of a manifold and its Lie algebra of vector 
fields have several natural constructions of modules coming from 
differential geometry. Examples of such modules include the space of 
functions on a manifold and more generally, the space of tensor fields.
In addition to these, one can also consider the spaces of jets of tensor
fields which also admit a natural action of the group of diffeomorphisms
and the Lie algebra of vector fields [S].

In this section we will review the construction of the bundle of jets of 
tensor fields and the module structure on the vector space of its sections.

First, let us discuss the notations that will be used in the paper.
We denote the set of non-negative integers by $\Z_+$, and consider a 
partial order on $\Z_+^n$ where $\alpha \geq \beta$ whenever $\alpha - 
\beta \in \Z_+^n$. The standard basis of $\Z^n$ will be denoted by
$\{\epsilon_1,\ldots,\epsilon_n \}$.
For a coordinate system $\{ x^i \}_{i=1,\ldots,n}$ and 
$\alpha = (\alpha_1, \ldots, \alpha_n) \in \Z_+^n$, we denote by 
$f^{(\alpha)}(x)$ the partial derivative 
$\left({\d \over \d x^1} \right)^{\alpha_1} \ldots 
\left({\d \over \d x^n} \right)^{\alpha_n} f(x)$. 
We will use notations $|\alpha| = \alpha_1 + \ldots + \alpha_n$, \ 
$\alpha ! = \alpha_1 ! \ldots \alpha_n !$, \ 
$\pmatrix{\alpha \cr \beta} = {\alpha! \over \beta! (\alpha - \beta)!}$, \
$(x-p)^\alpha = (x^1 - p^1)^{\alpha_1} \ldots (x^n - p^n)^{\alpha_n}$, etc.
We will also use a convention of dropping the summation
symbol when we take a sum over matching upper and lower indices:
$ u^i {\d \over \d x^i} = \sum\limits_{i=1}^n u^i {\d \over \d x^i}.$

 Let $M$ be an $n$-dimensional $C^\infty$ manifold, which is allowed to be
either real or complex. However the spaces of functions on $M$, vector fields, etc., 
will be always taken to be complex spaces. 

We begin by recalling the definition of an $N$-jet of a function. 
Let $f_1, f_2$ be two $C^\infty$ functions defined in a neighbourhood of 
a point $p$ on manifold $M$.
We say that $f_1$ is equivalent to $f_2$ at $p$, $f_1 \sim f_2$, if all
partial derivatives of $f_1$ and $f_2$ at $p$ of orders up to $N$ are
equal:
$$f_1^{(\alpha)} (p) = f_2^{(\alpha)} (p), 
\hbox{\rm {\hskip 1cm} for all \ } 0 \leq |\alpha| \leq N.$$
An equivalence class for this relation is called an 
$N$-jet of a function at $p$. 
Any $N$-jet at $p$ has a unique Taylor polynomial representative:
$$ \sum\limits_{0 \leq |\alpha| \leq N}
f_{[\alpha]} (p) (x-p)^\alpha .$$
The set of $N$-jets at $p$ forms a finite-dimensional vector
space with coordinates 
\break
$\left\{ f_{[\alpha]} (p) \; \big| \; 0 \leq |\alpha| \leq N \right\}$
and a basis 
$\left\{ (x-p)^\alpha \; \big| \; 0 \leq |\alpha| \leq N \right\}$.
If we now let $p$ vary over $M$, we get the vector bundle of $N$-jets of 
functions. Let us denote by $J_N (M)$ the space of sections of this bundle.

Now we are going to describe the action of the group of diffeomorphisms
$\Diff$ on $J_N (M)$. Suppose for a diffeomorphism $\varphi \in \Diff$ we
have $\varphi (p) = q$. Set $\psi = \varphi^{-1}$. Let $\{ x^i \}$ be 
a coordinate system near $p$, and  $\{ y^i \}$ be a coordinate system near $q$.

Let $F$ be a section of an $N$-jet bundle, $F \in J_N (M)$, with its value 
at $p$ given by the equivalence class
$$F(p) \sim \sum\limits_{0 \leq |\alpha| \leq N}
f_{[\alpha]} (p) (x-p)^\alpha .$$
Then the value of the section $\varphi F$ at point $q$ is defined by the jet
$$\varphi F(q) \sim \sum\limits_{0 \leq |\alpha| \leq N}
f_{[\alpha]} (p) (\psi(y)-\psi(q))^\alpha $$
$$\hbox{\hskip 1 cm} \sim \sum\limits_{0 \leq |\alpha| \leq N}
f_{[\alpha]} (p) 
\left( \sum\limits_{\beta > 0} {1 \over \beta !} 
\psi^{(\beta)} (q) (y-q)^\beta \right)^\alpha . \eqno{(1.1)}$$
The coordinates of $\varphi F$ at $q$ are then computed by expanding the last
expression in powers of $(y-q)$, and dropping terms of degrees greater than $N$.

If we pass to the infinitesimal action, we will obtain the action of the Lie 
algebra of vector fields $\Vect$ on $J_N (M)$.

For a vector field $\u = u^i(x) {\d \over \d x^i}$ we consider the corresponding
flow $\varphi_\epsilon \in \Diff$,
$$\varphi_\epsilon^i = x^i - \epsilon u^i (x) + o (\epsilon).$$
We now assume that the point $p$ and its image $q$ under $\varphi_\epsilon$
are in the same chart with coordinates $\{ x^i \}$.
Then for the inverse, $\psi_\epsilon = \varphi^{-1}_\epsilon \in \Diff$, we have
$$\psi_\epsilon^i = x^i + \epsilon u^i (x) + o (\epsilon),$$ 
and (1.1) becomes
$$\varphi_\epsilon F(q) \sim \sum\limits_{0 \leq |\alpha| \leq N}
f_{[\alpha]} (\psi_\epsilon (q)) 
\left( \sum\limits_{\beta > 0} {1 \over \beta !} 
\psi_\epsilon^{(\beta)} (q) (x-q)^\beta \right)^\alpha $$
$$\hbox{\hskip 1 cm} \sim F(q) 
+ \epsilon  \sum\limits_{0 \leq |\alpha| \leq N}
u^i(q) {\d f_{[\alpha]} \over \d x^i} (q) (x-q)^\alpha$$
$$+ \epsilon  \sum\limits_{0 \leq |\alpha| \leq N}
 f_{[\alpha]} (q) \sum\limits_{j=1}^n \alpha_j
\sum\limits_{\beta > 0} {1 \over \beta !} (u^j)^{(\beta)}(q)
(x-q)^{\alpha-\epsilon_j +\beta} + o(\epsilon) .$$
 Thus the action of $\Vect$ on $J_N(M)$ is given by the formula
$$\u F (q) \sim \sum\limits_{0 \leq |\alpha| \leq N}
 u^i(q) {\d f_{[\alpha]} \over \d x^i} (q) (x-q)^\alpha$$
$$+ \sum\limits_{0 \leq |\alpha| \leq N}
 f_{[\alpha]} (q) \sum\limits_{j=1}^n \alpha_j
\sum\limits_{\beta > 0} {1 \over \beta !} (u^j)^{(\beta)}(q)
(x-q)^{\alpha-\epsilon_j +\beta}. \eqno{(1.2)} $$  

 We can see from these formulas that the subspace of sections with derivatives
up to order $\ell$ vanishing everywhere
$$ \left\{ F \in J_N(M) \; \big| \; \; f_{[\alpha]} \equiv 0 
\hbox{\rm {\hskip 0.3cm} for all \ } 0 \leq |\alpha| \leq \ell \right\}$$
is a submodule for the actions of both $\Diff$ and $\Vect$. The factor-module
in this case is isomorphic to the space $J_\ell (M)$ of sections of the bundle 
of $\ell$-jets.

 Now let us look at a more general case of tensor fields. A tensor
field of type $(s,k)$ is a section of the corresponding tensor bundle. Suppose
that a tensor field $F$ is given in a chart with local coordinates $\{ x^i \}$
by expression
$$ F(x) = f^{(i_1 \ldots i_s)}_{(j_1 \ldots j_k)} (x) dx^{j_1} \ldots dx^{j_k}
{\d \over \d x^{i_1}} \ldots {\d \over \d x^{i_s}} .$$
Then the diffeomorphism $\varphi \in \Diff$ with $\varphi^{-1} = \psi$ acts
on $F$ according to the formula
$$ \varphi F (y) = \hfill $$
$$ f^{(i_1 \ldots i_s)}_{(j_1 \ldots j_k)} (\psi(y))
{\d \psi^{j_1} \over \d y^{j_1^\prime}} \ldots {\d \psi^{j_k} \over \d y^{j_k^\prime}} 
{\d \varphi^{i_1^\prime} \over \d x^{i_1}} \left( \psi(y) \right) \ldots 
{\d \varphi^{i_s^\prime} \over \d x^{i_s}} \left( \psi(y) \right)
d y^{j_1^\prime} \ldots d y^{j_k^\prime} {\d \over \d y^{i_1^\prime}} \ldots  
{\d \over \d y^{i_s^\prime}} . \eqno{(1.3)}$$
The corresponding action of the Lie algebra $\Vect$ can be conveniently encoded
using representations of the Lie algebra $\gl$. 
Let us explain this construction. The tensor bundle in question is a 
tensor product of $s$ copies of the tangent bundle and $k$ copies of the cotangent
bundle. For a given coordinate system $\{ x^i \}$ there is an action of the Lie 
algebra $\gl$ on the cotangent space, where we set $\{d x^i\}$
as the standard basis of the natural $\gl$-module. The tangent space becomes
the conatural module for $\gl$. The action of the elementary matrices $E^p_q$
(a matrix that has entry $1$ in position $(p,q)$ and zeros elsewhere) 
on the tangent and cotangent spaces is given by the formulas:
$$ E^p_q dx^i = \delta_q^i dx^p ,$$
$$ E^p_q {\d \over \d x^i} = - \delta^p_i {\d \over \d x^q}. $$
The fiber $V = \Span \left< dx^{j_1} \ldots dx^{j_k} {\d \over \d x^{i_1}} \ldots
{\d \over \d x^{i_s}} \big| 1 \leq j_1, \ldots, j_k, i_1, \ldots i_s \leq n \right>$
of the $(s,k)$ tensor bundle gets the structure of a $\gl$-module as a tensor product.

In local coordinates $\{ x^i \}$, any tensor field of type $(s,k)$ is a linear 
combination of $f(x) v$, where $v \in V$. The action of a vector field 
$\u = u^i(x) {\d \over \d x^i}$ is then given by the (Lie derivative) formula:
$$\u \left( f(x) v \right) = \left( u^i(x) {\d f \over \d x^i} \right) v
+ \left( f(x) {\d u^i \over \d x^j} \right) E_i^j v .  \eqno{(1.4)} $$   
Just as in the case of functions, we can introduce the bundle of $N$-jets of
tensor fields. The fiber at a point $p$ is spanned by the jets
$$ v^{(\alpha)} = (x-p)^\alpha v, \; \; v\in V, \; 0\leq |\alpha| \leq N.$$
Again, let us define the action of $\Diff$ on the space $J_N T^{(s,k)} (M)$ of
sections of this bundle. If we have a section $F$ with value at a point $p$
given by an equivalence class of a local tensor field
$$ F(p) \sim \sum\limits_{0\leq |\alpha| \leq N} 
f_{(j_1\ldots j_k)[\alpha]}^{(i_1\ldots i_s)} (p) (x-p)^\alpha
dx^{j_1} \ldots dx^{j_k} {\d \over \d x^{i_1}} \ldots {\d \over \d x^{i_s}}$$
then $\varphi F$ is a section with value at $q = \varphi(p)$ given by the equivalence
class of the image of this local tensor field under the action of $\varphi$
according to (1.3):
$$ \varphi F(q) \sim 
\sum\limits_{0\leq |\alpha| \leq N} 
f_{(j_1\ldots j_k)[\alpha]}^{(i_1\ldots i_s)} (\psi(q)) (\psi(y)-\psi(q))^\alpha \times$$
$$\times
{\d \psi^{j_1} \over \d y^{j_1^\prime}} \ldots {\d \psi^{j_k} \over \d y^{j_k^\prime}} 
{\d \varphi^{i_1^\prime} \over \d x^{i_1}} \left( \psi(y) \right) \ldots 
{\d \varphi^{i_s^\prime} \over \d x^{i_s}} \left( \psi(y) \right)
d y^{j_1^\prime} \ldots d y^{j_k^\prime} {\d \over \d y^{i_1^\prime}} \ldots  
{\d \over \d y^{i_s^\prime}} .$$
Infinitesimal variant of the above formula yields the action of the Lie algebra
$\Vect$ on $J_N T^{(s,k)} (M)$:
$$\u \left( f(x) v^{(\alpha)} \right) = 
\left( u^i(x) {\d f(x) \over \d x_i} \right) v^{(\alpha)}  
+ \sum\limits_{j=1}^n \alpha_j \sum\limits_{\beta > 0} {1 \over \beta !}
f(x) (u^j)^{(\beta)} (x) v^{(\alpha - \epsilon_j + \beta)}$$
$$ + \sum\limits_{j=1}^n \sum\limits_{\beta \geq 0} {1 \over \beta !}
f(x) (u^i)^{(\beta + \epsilon_j)} (x) 
\left( E_i^j v \right)^{(\alpha + \beta)}. \eqno{(1.5)} $$

 We would like to point out some properties of the module $J_N T^{(s,k)} (M)$.
In addition to being a module for the Lie algebra $\Vect$, it is also a module
over a commutative algebra $\F (M)$ of functions on $M$. Moreover (1.5) shows
that the two structures are compatible in the following way:
$$\u \left( f(x) F \right) = \left( \u f(x) \right) F + f(x) \left( \u F \right),
\; \; \u\in\Vect, f\in \F(M), F \in J_N T^{(s,k)} (M) .$$
We see that vector fields act as derivations of the multiplication of the jets 
of tensor fields by functions. This motivates the definition of a category 
of modules that will be introduced in the next section.

\

\

{\bf 2. Category $\J$}.

\

For the rest of the paper the manifold will be an $n$-dimensional torus,
$M = \T = \R^n / \Z^n$. The algebra of functions on a torus has the Fourier basis 
$\left\{ \e{m} \big| m\in\Z^n \right\}$. The Lie algebra $\VT$, also denoted
as $W_n$, is a free module over $\FT$ of rank $n$ with a basis (as an 
$\FT$-module) $\left\{ d_j = {1 \over 2\pi i} {\d \over \d x^j} \big| 
j=1,\ldots,n \right\}$. A basis of $W_n$ over $\C$ is given by
$$\left\{ d_j (s) = {1 \over 2\pi i} \e{s} {\d \over \d x^j} \big| 
s\in\Z^n, j=1,\ldots,n \right\} .$$
The subspace spanned by $\{ d_j \}_{j=1,\ldots,n}$ is a Cartan subalgebra
in $W_n$.

Note that because the torus is a flat manifold, 
all bundles considered in Section 1 are trivial.
Equivalently, the spaces of sections of these bundles are free modules
over $\FT$.

Let us define the following category $\J$ of $W_n$-modules:

{\bf Definition.}
A $W_n$-module $J$ belongs to category $\J$ if the following properties hold:

\noindent
(J1) The action of $d_j$, $j=1,\ldots,n,$ on $J$ is diagonalizable.

\noindent
(J2) $J$ is a free $\FT$-module of a finite rank.

\noindent
(J3) For any $\u\in W_n, f\in\FT, w\in J,$
$$ \u \left( f w \right) = \left( \u f \right) w + f \left( \u w \right) .$$

 Strictly speaking, we should denote this category $\J_n$, but we will omit
the subscript most of the time. All the $\VT$-modules discussed in Section 1
belong to category $\J$.

The goal of this paper is to classify indecomposable and irreducible modules
in this category. Irreducible modules in $\J$ have been already classified by
Eswara Rao [R2], but here we simplify the proof and make it more lucid.
A class of indecomposable $W_n$-modules has been constructed in [R1].

{\bf Remark.} When we talk about a submodule for a module in category $\J$
we mean a subspace which is invariant under the action of both $W_n$ and $\FT$.
Also, in this paper an indecomposable module will be understood to be non-zero. 

The concept of a polynomial module will play a central role in our proof. The
definition of a polynomial module was given in [BB] in a general setup. Here
we adapt that definition for the particular case of $W_n$. We will show below
that one can choose a basis $v_1, \ldots v_k$ of $J$ over $\FT$ in
such a way that the action of $W_n$ is as follows:
$$ d_j (s) \left( \e{m} v_r \right) = \sum\limits_{\ell=1}^k f_{jr\ell} (s,m)
\e{(m+s)} v_\ell , \; \; s,m\in\Z^n . \eqno{(2.1)}$$
 
We say that module $J$ is {\it a polynomial module} if the structure constants
$f_{jr\ell} (s,m)$ are polynomials in $s,m\in\Z^n$.

The desired classification of modules in $\J$ will be obtained in three steps
by proving the following theorems (see the next section for more detailed
statements).

\

{\bf Theorem 1.} Let $n=1$. Every $W_1$-module in category $\J_1$ is a 
polynomial module.

 We will see that with very little effort one can derive from Theorem 1 its 
generalization to an arbitrary rank:

{\bf Theorem 2.} Every $W_n$-module in category $\J$ is a polynomial module.

From this Theorem we almost immediately deduce the classification of the
indecomposable modules in category $\J$. In order to state this result, we
would need to introduce the Lie algebra $W_n^+$.

 Consider the Lie algebra of derivations of the algebra of polynomials
in $n$ variables:
$$\Der \left( \C [z_1, \ldots, z_n] \right)
= \Span \left< z^\alpha {\d \over \d z_j} \; \bigg| \; \alpha\in\Z_+^n, 
j=1,\ldots,n \right> .$$
The Lie algebra $W_n^+$ is defined as a subalgebra in 
$\Der \left( \C [z_1, \ldots, z_n] \right)$:
$$ W_n^+ = \Span \left< z^\alpha {\d \over \d z_j} \; \bigg| \; 
\alpha\in\Z_+^n \backslash \{ 0 \}, j=1,\ldots,n \right> .$$

{\bf Theorem 3.} There is a 1-1 correspondence between indecomposable modules
in category $\J$ and pairs $(\lambda, U)$, where $\lambda \in \C^n / \Z^n$
and $U$ is an indecomposable finite-dimensional module for $W_n^+$.

\

\

{\bf 3. Structure of indecomposable modules.}

\

 Finiteness condition (J2) implies that every module in category $\J$ is a
finite direct sum of indecomposable submodules. This allows us to restrict
our attention to indecomposable modules.

 Let us write the weight decomposition of an indecomposable module $J$
with respect to the Cartan subalgebra of $W_n$:
$$J = \mathop\oplus\limits_{\mu\in\C^n} J_\mu, \; \; {\rm where} \;
J_\mu = \left\{ w\in J \big| d_j w = \mu_j w \right\} .$$
It is easy to see that $d_j (s) \J_\mu \subseteq J_{\mu + s}$
and $\e{s} \J_\mu \subseteq J_{\mu + s}$, thus the weights of $J$ are 
split
into $\Z^n$-cosets in $\C^n$, and the submodules corresponding to distinct
cosets form a direct sum. Since we assumed $J$ to be indecomposable, its
weight lattice is a single coset $\lambda + \Z^n, \lambda\in\C^n$. 
Let us denote the weight space $J_\lambda$ by $U$. It is easy to see that
the basis of $U$ is also a basis of $J$ as an $\FT$-module. Thus, by (J2),
the space $U$ is finite-dimensional and 
$ J = \FT \otimes U .$

Since $\e{s} U = J_{\lambda + s}$, we may identify each weight space of $J$
with the same finite-dimensional space $U$. The operator
$$ d_j(s) : U \rightarrow \e{s} U $$
induces an endomorphism $D_j(s) : \; U \rightarrow U$, such that
$$ d_j(s) \big|_U = \e{s} D_j(s) .$$
In particular we have $D_j(0) = \lambda_j \Id$. 

The finite-dimensional operator $D_j(s)\in \End U$ completely determines the action
of $d_j(s)$ on $J$ since by (J3),
$$d_j(s) \left( \e{m} v \right) = \left( d_j(s) \e{m} \right) v 
+ \e{m} d_j(s) v = 
\e{(m+s)} \left( m_j \Id + D_j(s) \right) v . \eqno{(3.1)}$$

% Let us summarize our discussion in the following preliminary
% 
% {\bf Lemma 1.} Every indecomposable module $J$ in category $\J$ can be
% written as a tensor product
% $$ J = \FT \otimes U , \; \; \; \dim U < \infty, $$
% where the action of $W_n$ is given by 
% $$d_j(s) \left( \e{m} v \right) = \e{(m+s)} \left( m_j \Id + D_j(s) \right) v,
% \eqno{(3.1)}$$
% for some $D_j (s) \in \End U$ with $D_j (0) = \lambda_j \Id$ .

Now we see that the action of $d_j (s)$ is written in the form (2.1)
The structure constants $f_{jr\ell}$ are encoded in the operators
$m_j \Id + D_j (s)$. The dependence on $m$ here is clearly polynomial, so
the statement that $J$ is a polynomial $W_n$-module is equivalent to the claim
that the dependence of the family of operators $\{ D_j (s) \}$ on $s$ is 
polynomial.

\

{\bf Theorem 2.} An indecomposable $W_n$-module $J = \FT \otimes U$ in 
category $\J$ is a polynomial module. The action of $W_n$ can be written
as follows:
$$d_j(s) \left( \e{m} v \right) = \e{(m+s)} \left( m_j \Id + D_j(s) \right) v,$$
where operators $D_j(s) \in \End U$ have a polynomial dependence on $s\in\Z^n$
and $D_j (0) = \lambda_j \Id$.

Proving Theorem 2 directly would be rather technical. 
% (see [R2] where the case of irreducible modules is treated)
Instead, we will derive it from its special
case of $n=1$, which is Theorem 1. The proof of Theorem 1 will be deferred to
the next section.

First of all, let us write down the commutator relations between operators
$D_j(s)$:

{\bf Lemma 1.} 
$$ \left[ D_j (s), D_k (m) \right] = m_j \left( D_k (s+m) - D_k (m) \right)
- s_k \left( D_j (s+m) - D_j(s) \right) . \eqno{(3.2)} $$

This Lemma can be derived in a straightforward way from the commutator
relations in $W_n$
$$  \left[ d_j (s), d_k (m) \right] = m_j d_k (s+m) - s_k d_j (s+m) \eqno{(3.3)}$$
and (3.1). The details of this computation are left as an exercise.

{\sl Proof of Theorem 2.}
 We assume that the claim of the theorem holds in rank one case, $n=1$.

Consider the following subalgebras in $W_n$, each isomorphic to $W_1$:
$$ W_1^{(j)} = \Span \left< {1\over 2\pi i} e^{2 \pi i s x^j} {\d \over \d x^j} 
\; \bigg| \; s\in\Z \right> .$$
The subspace 
$ \mathop\oplus\limits_{m\in\Z} e^{2\pi i m x^j} U \subset J$
is a $W_1^{(j)}$-module which belongs to category $\J_1$. Applying Theorem 1,
we get that the family of operators $\{ D_j (s \epsilon_j) \}$ has a polynomial
dependence on $s\in \Z$.

Without the loss of generality we may restrict ourselves to proving that
$D_1(s_1, \ldots, s_n)$ is a polynomial in $s_1,\ldots,s_n$. This is 
of course equivalent to showing that 
\break
$D_1(s_1, 1+s_2, \ldots, 1+s_n)$ is a 
polynomial in $s_1, s_2,\ldots,s_n$. We will prove by induction the claim
that $D_1(s_1, 1+s_2, \ldots, 1+s_j, 1, \ldots, 1)$ is a 
polynomial in $s_1, s_2,\ldots,s_j$.

Let us establish the basis of induction. From (3.2) we get that
$$ \left[ D_1 (s_1 \epsilon_1), D_1 (0, 1, \ldots 1) \right] = 
-s_1 D_1(s_1, 1, \ldots ,1) + s_1 D_1(s_1 \epsilon_1) ,$$
and so 
$$ s_1 D_1(s_1, 1, \ldots ,1) = s_1 D_1(s_1 \epsilon_1)
- \left[ D_1 (s_1 \epsilon_1), D_1 (0, 1, \ldots 1) \right] .$$
The right hand side is manifestly a polynomial in $s_1$, so is the left
hand side. Note also that the right hand side vanishes at $s_1 = 0$
because $D_1(0) = \lambda_1 \Id$. Hence this polynomial has a factor of $s_1$,
and thus $D_1(s_1, 1, \ldots ,1)$ is a polynomial in $s_1$, which proves
the basis of induction. 

 Let us now prove the inductive step. Again from (3.2) we get
% $$ \left[ D_j (s_j \epsilon_j), 
% D_1 (s_1, 1 + s_2, \ldots 1+s_{j-1}, 1, \ldots, 1) \right] =
% D_1 (s_1, 1 + s_2, \ldots 1+s_{j-1}, 1+s_j, 1, \ldots, 1)
% - D_1 (s_1, 1 + s_2, \ldots 1+s_{j-1}, 1, \ldots, 1) ,$$
% and so 
$$ D_1 (s_1, 1 + s_2, \ldots 1+s_{j-1}, 1+s_j, 1, \ldots, 1) =$$
$$ \left[ D_j (s_j \epsilon_j), 
D_1 (s_1, 1 + s_2, \ldots 1+s_{j-1}, 1, \ldots, 1) \right] 
+ D_1 (s_1, 1 + s_2, \ldots 1+s_{j-1}, 1, \ldots, 1) .$$
By induction assumption, the right hand side is a polynomial in
$s_1, \ldots, s_{j-1}, s_j$, and so is the left hand side. This 
completes the induction and Theorem 2 is now proved (under assumption of
validity of Theorem 1).

Now given that $D_j (s)$ is a polynomial function with a constant term
$D_j (0) = \lambda_j \Id$, we can expand it in a finite sum:
$$ D_j (s) = \lambda_j \Id + 
\sum\limits_{\alpha\in\Z_+^n \backslash \{ 0 \}}^{\rm finite}
{s^\alpha \over \alpha!} D_j^{(\alpha)} , \eqno{(3.4)}$$
where the operators $D_j^{(\alpha)} \in\End U$ are independent of $s$
and $D_j^{(\alpha)} = 0$ for all but finitely many 
$\alpha\in\Z_+^n \backslash \{ 0 \}$.

\

{\bf Theorem 3.} 
(a) The operators $D_j^{(\alpha)} \in\End U$ yield a
finite-dimensional representation on space $U$ of a Lie algebra
$$ W_n^+ = \Span \left< z^\alpha {\d \over \d z_j} \; \bigg| \; \; 
\alpha\in\Z_+^n \backslash \{ 0 \}, \; j=1,\ldots,n \right> ,$$
given by $\rho \left( z^\alpha {\d \over \d z_j} \right) = D_j^{(\alpha)}$.

(b) There is a 1-1 correspondence between indecomposable modules in 
category $\J$ and pairs $(\lambda, U)$, where $\lambda\in \C^n / \Z^n$
and $U$ is an indecomposable finite-dimensional module for $W_n^+$.
This correspondence is given by the tensor product decomposition
$J = \FT \otimes U$, and the action of $W_n$ is
$$ {1 \over 2 \pi i} \e{s} {\d \over \d x^j} 
\left( \e{m} v \right) 
= (m_j + \lambda_j) \e{(m+s)} v 
+ \sum\limits_{\beta > 0} {s^\beta \over \beta !} \e{(m+s)} 
\rho \left( z^\beta {\d \over \d z_j} \right) v .\eqno{(3.5)}$$

(c) There is a 1-1 correspondence between irreducible modules in 
category $\J$ and pairs $(\lambda, V)$, where $\lambda\in \C^n / \Z^n$
and $V$ is an irreducible $\gl (\C)$-module. The action of $W_n$ on
$J = \FT \otimes V$ is as follows:
$$ {1 \over 2 \pi i} \e{s} {\d \over \d x^j} 
\left( \e{m} v \right) 
= (m_j + \lambda_j) \e{(m+s)} v 
+ \sum\limits_{p=1}^n s_p \e{(m+s)} 
\rho \left( E^p_j \right) v .\eqno{(3.6)}$$

{\sl Proof.} Let us determine the commutator relations between operators
$D_j^{(\alpha)}$. We have 
$$ \left[ D_j (s), D_k (m) \right] =
\sum\limits_{\alpha,\beta\in\Z_+^n\backslash\{ 0 \}}
{s^\alpha m^\beta \over \alpha ! \beta !} 
\left[ D_j^{(\alpha)}, D_k^{(\beta)} \right] .$$
On the other hand, by Lemma 1,
$$ \left[ D_j (s), D_k (m) \right] =
\sum\limits_{\gamma\in\Z_+^n\backslash\{ 0 \}}
m_j {(s+m)^\gamma - m^\gamma \over \gamma !} D_k^{(\gamma)}
- \sum\limits_{\gamma\in\Z_+^n\backslash\{ 0 \}}
s_k {(s+m)^\gamma - s^\gamma \over \gamma !} D_j^{(\gamma)} .$$

Two polynomials have equal values whenever their coefficients coincide.
If we equate the coefficients at ${s^\alpha m^\beta \over \alpha ! \beta !},
\alpha, \beta \in \Z_+^n\backslash\{ 0 \}$, we get
$$ \left[ D_j^{(\alpha)}, D_k^{(\beta)} \right] =
\beta_j D_k^{(\alpha+\beta-\epsilon_j)} 
- \alpha_k D_j^{(\alpha+\beta-\epsilon_k)} . \eqno{(3.7)}$$
We see that these are precisely the commutator relations in the algebra
$W_n^+$:
$$ \left[ z^\alpha {\d \over \d z_j}, z^\beta {\d \over \d z_k} \right]
= \beta_j z^{\alpha+\beta-\epsilon_j} {\d \over \d z_k}
- \alpha_k z^{\alpha+\beta-\epsilon_k} {\d \over \d z_j} ,$$
so the map $\rho : \; W_n^+ \rightarrow \End U$, given by
$\rho\left( z^\alpha {\d \over \d z_j} \right) = D_j^{(\alpha)}$,
is a representation of $W_n^+$. This completes the proof of part (a).

 Let us prove (b). We have already seen that an indecomposable module
$J$ in category $\J$ yields a coset of weights $\lambda + \Z^n \subset \C^n$,
and a finite-dimensional representation of $W_n^+$. It is easy to check that 
the $W_n^+$-module $U$ is independent of the choice of the weight $\lambda$ 
in the coset. 
 Conversely, the commutator relations (3.7) imply (3.2) and together with (3.1)
give (3.3), provided of course that the right hand side in (3.4) is finite.
Thus we need to show that for a finite-dimensional representation $(U,\rho)$
of the Lie algebra $W_n^+$, we will have that 
$D_j^{(\alpha)} = \rho\left( z^\alpha {\d \over \d z_j} \right) = 0$ for all 
but finitely many $\alpha\in\Z_+^n  \backslash\{ 0 \}$. This will follow from
the following simple

{\bf Lemma 2.} Let $\L$ be a Lie algebra, and $(U,\rho)$
its finite-dimensional representation. Suppose that for $x, y_1, y_2, \ldots
\in\L$ we have
$$[x, y_k ] = \nu_k y_k, \; \; \nu_k\in\C \; k=1,2,\ldots . $$
Then there are at most $\left( \dim U \right)^2 - \dim U + 1$
distinct eigenvalues for which $\rho(y_k) \neq 0$.

{\sl Proof of the Lemma.} In representation $\rho$ we have
$$[\rho(x), \rho(y_k)] = \nu_k \rho(y_k) .$$
However an element $\rho(x)$ in the Lie algebra $gl(U)$ may have at most
$\left( \dim U \right)^2 - \dim U + 1$ distinct eigenvalues in the adjoint
representation, which implies the claim of the Lemma.

 To apply this lemma, we consider the element $E = z_1 {\d \over \d z_1} +
\ldots z_n {\d \over \d z_n} \in W_n^+$ . The following relations hold:
$$ \left[ E, z^\alpha {\d \over \d z_j} \right] =
\left( |\alpha| - 1\right) z^\alpha {\d \over \d z_j} .$$
Then by Lemma 2, $\rho\left( z^\alpha {\d \over \d z_j} \right) = 0$
for all but finitely many $\alpha$.

 To complete the proof of part (b) of the theorem, we note that an 
indecomposable $W_n$-module in category $\J$ yields an indecomposable
$W_n^+$-module and vice versa.

 Now let us prove part (c) and assume that $J$ is an irreducible module
in category $\J$. From part (b) we know that such a module is determined 
by a pair $(\lambda, V)$, where $\lambda \in \C^n / \Z^n$ and $V$ is a 
finite-dimensional $W_n^+$-module. The module $V$ has to be irreducible,
because we have a correspondence between $W_n^+$-submodules $S \subset V$
and $W_n$-submodules $\FT \otimes S \subset J$.

Let us show that irreducible finite-dimensional $W_n^+$-modules $V$ are
just irreducible $\gl(\C)$-modules. Let $V_\nu$ be an eigenspace of 
the operator $\rho(E)$ in $V$ corresponding to eigenvalue $\nu$. It is easy to
see that
$$\rho\left(z^\alpha {\d \over \d z_j} \right) V_\nu \subseteq 
V_{\nu + |\alpha| -1} .$$
This implies that $\mathop\oplus\limits_{k=1}^\infty V_{\nu+k}$ is
a $W_n^+$-submodule in $V$. However $V$ is irreducible, which implies
that $V_{\nu+k} = (0)$ for all $k=1, 2, \ldots$, and hence
$$\rho\left(z^\alpha {\d \over \d z_j} \right) = 0 \; \; 
\hbox{\rm for all} \; |\alpha| > 1 .$$
Thus the ideal 
$$ W_n^{++} = \Span \left< z^\alpha {\d \over \d z_j} \; \bigg| \; \; 
\alpha\in\Z_+^n, \; |\alpha| > 1, \; j=1,\ldots,n \right> $$
vanishes in every finite-dimensional irreducible $W_n^+$-module.
But $W_n^+ / W_n^{++} \cong \gl(\C)$, and the claim of part (c) of the theorem
follows.

\

 Next, let us give an example of a family of indecomposable finite-dimensional
$W_n^+$-modules.

{\bf Example.} Let $V$ be an irreducible finite-dimensional $\gl(\C)$-module
and let $$\widetilde V = \C[z_1,\ldots,z_n] \otimes V .$$
We can define on $\widetilde V$ the structure of a tensor module for the Lie
algebra $\Der \C[z_1,\ldots,z_n]$ and its subalgebra 
$W_n^+$ (cf., (1.4), see also [Ru]):
$$ z^\beta {\d \over \d z_j} (z^\alpha v) = 
\alpha_j z^{\alpha+\beta-\epsilon_j} v
+ \sum_{k=1}^n \beta_k z^{\alpha+\beta-\epsilon_k} (E_j^k v). $$
As a $W_n^+$-module, this tensor module has finite-dimensional factors
$$ V^{(N)} = \widetilde V / \left< z^\alpha \otimes V \; \big| \; |\alpha| > N \right> .$$
Using the correspondence of Theorem 3(b), we construct a representation of
the Lie algebra $\VT$ on space $\FT \otimes V^{(N)}$ (setting $\lambda = 0$):
$${1\over 2\pi i} \e{s} {\d \over \d x^j} \left( \e{m} z^\alpha v \right)
= m_j \e{(m+s)} z^\alpha v $$ 
$$ + \alpha_j \sum\limits_{\beta > 0} 
{s^\beta \over \beta!} \e{(m+s)} z^{\alpha+\beta-\epsilon_j} v 
+ \sum\limits_{\beta > 0} {s^\beta \over \beta!}
\sum\limits_{k=1}^n \beta_k \e{(m+s)} z^{\alpha+\beta-\epsilon_k} E_j^k v .
\eqno{(3.8)} $$
% Here we denoted $v^{(\alpha)} = (2\pi i)^{|\alpha|} z^\alpha v$. 

Comparing (3.8) with (1.5) we see that the module $\FT \otimes V^{(N)}$ is in 
fact isomorphic to a module of sections of $N$-jets of tensor fields
corresponding to the $\gl$-module $V$.
The isomorphism is given by $z^\alpha v = (2 \pi i)^{|\alpha|} v^{(\alpha)}$.

{\bf Remark.} It is curious to point out that for this family of polynomial
modules, the degree of the structure polynomials is equal to $N+1$, and thus 
could be arbitrarily high. In the previously known examples (see [BB] and
[BZ]) the degree was at most 3.

\

\

{\bf 4. Rank one case.}

\

In this section we will give the proof of Theorem 1. Since we will deal
exclusively with the case $n=1$, we may simplify notations denoting
$x^1$ as $x$, $d_1(s)$ as $d(s)$, etc. Let $J$ be an indecomposable 
module in category $\J_1$. As in Section 3 we see that 
$J = \F ({\BBB T}^1) \otimes U$. Then (3.1) becomes
$$d(s) \left( \e{m} v \right) = \e{(m+s)} \left( m \Id + D(s) \right) v ,
\; \; s,m\in\Z.  \eqno{(4.1)}$$
The operators $D(s)\in \End U$ satisfy the commutator relations (cf. (3.2)):
$$ \left[ D(s), D(m) \right] = (m-s) D(s+m) - m D(m) + s D(s) . \eqno{(4.2)}$$

{\bf Theorem 1.} Every $W_1$-module $J$ in category $\J_1$ is a polynomial module.  

{\sl Proof.} Just as in Section 3, we will assume $J$ to be indecomposable,
and so (4.1) holds. We will prove this theorem by showing that a family of
operators $\{ D(s) \}$ on a finite-dimensional space $U$ satisfying 
(4.2) must have a polynomial dependence on $s\in\Z$.

The relations (4.2) define an infinite-dimensional Lie algebra $\L$ with
basis 
\break
$\left\{ D(s) \big| s\in\Z \right\}$, and we are studying a 
finite-dimensional representation $\rho$ of $\L$ on space $U$.

Our strategy will be the following. First we will find eigenvectors of $D(-1)$
in the adjoint representation of $\L$. We shall see that these eigenvectors
are difference derivatives of $D(s)$ with respect to $s$. By applying Lemma 2,
we will conclude that higher order derivatives of $D(s)$ vanish, which means 
that $D(s)$ is a polynomial in $s$.

Let us define the operator of a difference derivative. Let
$$ f: \; \Z \rightarrow A$$
be a function of an integer variable with values in an abelian group $A$
(in our case $A$ is the vector space $\L$). The {\it difference derivative}
$\d f$ is a function $\d f: \; \Z \rightarrow A$, defined by
$\d f (s) = f(s+1) - f(s)$. By iteration, we can also define higher order 
difference derivatives of $f$. It is easy to see that
$$\d^m f(s) = \sum\limits_{k=0}^m (-1)^{m-k} \pmatrix{m \cr k \cr} f(s+k) .
\eqno{(4.3)} $$
Let us now additionally assume that $A$ is a module over $\Q$, so that we can
define interpolation polynomials. For any $N+1$ distinct points $r_1,\ldots,
r_{N+1} \in\Z$ and any set of values $a_1, \ldots, a_{N+1} \in A$, there
exists a unique polynomial $f(t)\in A[t]$ of degree at most $N$ such
that $f(r_j) = a_j$ for $j=1,\ldots, N+1$ ([vdW], section 22).

\

{\bf Lemma 3.} Fix $s\in \Z$.

(a) If $\d^m f (s) = 0$ for all $m \geq 0$
then $f(r) = 0$ for all $r \geq s$.

(b) Suppose $\d^m f(s) = 0$ for all $m > N$. Let $g(t)$ be the interpolation
polynomial of degree at most $N$ defined by $g(j) = f(j)$ for 
$j = s, s+1, \ldots, s+N$. Then $f(r) = g(r)$ for all $r\geq s$.

{\sl Proof of the Lemma.} Part (a) can be easily proved by induction 
using (4.3). To prove part (b), we consider the function $h(r) = f(r) - g(r)$.
Since $h(s) = \ldots = h(s+N) = 0$, we get from (4.3) that $\d^m h (s) = 0$
for $m = 0, \ldots, N$. For $m > N$ we have $\d^m f (s) = 0$ by assumption
and $\d^m g (s) = 0$ since it is a polynomial of degree at most $N$.
Thus $\d^m h (s) = 0$ for all $m \geq 0$, and according to part (a) of the
Lemma, $f(r) = g(r)$ for all $r\geq s$. This completes the proof of the 
Lemma.

 Now consider the following elements in the Lie algebra $\L$:
$$ y_m = \d^{m+1} D(-1) = 
\sum\limits_{k=0}^{m+1} (-1)^{m+1-k} \pmatrix{m+1 \cr k \cr} D(-1+k) .$$

{\bf Lemma 4.} (a) For $m \geq 0$, $y_m$ is an eigenvector for $\ad D(-1)$
with eigenvalue 
\break
$\nu_m = -m$.

(b) For $m,k \geq 0$, \ $[y_k, y_m] = (m-k) y_{m+k}$.

{\sl Proof.} In the calculation below we will use two elementary properties
of binomial coefficients:
$$\sum\limits_{k=0}^m (-1)^{m-k} \pmatrix{ m \cr k \cr} = 0 \; \hbox{\rm for \ }
m\geq 1, 
\hbox{\rm \ \ and \ \ }
(r+1) \pmatrix{ m \cr r+1 \cr} = (m-r) \pmatrix{ m \cr r \cr}. $$

Now,
$$\left[ D(-1), y_{m-1} \right] =
\sum\limits_{k=0}^m (-1)^{m-k} \pmatrix{m \cr k \cr} [D(-1), D(-1+k)] $$
$$ = \sum\limits_{k=0}^m (-1)^{m-k} \pmatrix{m \cr k \cr} 
\bigg( k D(-2+k) - (k-1) D(-1+k) - D(-1) \bigg)$$
$$ = - \sum\limits_{{r=0} \atop {r=k-1}}^{m-1} (-1)^{m-r} (r+1) \pmatrix{m \cr r+1 \cr}
D(-1+r)
- \sum\limits_{k=0}^m (-1)^{m-k} (k-1) \pmatrix{m \cr k \cr} D(-1+k)$$
$$- \sum\limits_{k=0}^m (-1)^{m-k} \pmatrix{m \cr k \cr} D(-1) $$
$$ = - \sum\limits_{r=0}^m (-1)^{m-r} (m-r) \pmatrix{m \cr r \cr} D(-1+r)
- \sum\limits_{r=0}^m (-1)^{m-r} (r-1) \pmatrix{m \cr r \cr} D(-1+r) $$
$$ = (-m+1) \sum\limits_{r=0}^m (-1)^{m-r} \pmatrix{m \cr r \cr} D(-1+r)
 = (-m+1) y_{m-1} .$$

% To prove part (b) we will define the leading term of an element 
% $\sum\limits_{r=-1}^m a_r D(r)$ (with $a_m \neq 0$) to be $a_m D(m)$.
% It follows from (4.2) that given two such elements with leading terms
% $a_m D(m)$ and $b_k D(k)$, the leading term of their commutator will be
% $(k-m) a_m b_k D(m+k)$, provided that $m\neq k$.
%  
%  As a consequence of part (a) we see that the operator $\ad D(-1)$ is 
% diagonalizable on the subspace 
% $S_N = \Span\left< D(j) \big| -1\leq j \leq N \right>$
% with the eigenbasis $D(-1), y_0, \ldots, y_N$ with eigenvalues $0, 0, 1, \ldots N$.
% Suppose $m$ and $k$ are distinct and $m+k \leq N$. The leading term of $y_m$
% is $D(m)$, and the leading term of $y_k$ is $D(k)$. Let us assume that $m\neq k$.
% Then the leading term of $[y_m, y_k]$ is $(k-m) D(m+k)$. Thus $[y_m, y_k] \in S_N$.
% Moreover, 
% $$ [D(-1), [y_m, y_k]] = [[D(-1), y_m], y_k] + [y_m, [D(-1), y_k]] = (m+k) [y_m, y_k] .$$
% Since $y_{m+k}$ is the only eigenvector in $S_N$ with eigenvalue $(m+k)$,
% we conclude that $[y_m, y_k]$ is a multiple of $y_{m+k}$. Comparing their leading terms,
% we see that $[y_m, y_k] = (k-m) y_{m+k}$, which is the claim of part (b) of the lemma.

 Part (b) of the lemma will not be used in this paper and its proof is left as an 
exercise.

\

Now we are ready to complete the proof of Theorem 1.
Combining Lemma 4(a) with Lemma 2, we conclude that there exists $N$ 
such that $\rho(y_m) = \rho(\d^{m+1} D(-1)) = 0$ 
for all $m \geq N$. Then by Lemma 3(b), there exists an $\End U$-valued polynomial
$g(s)$ such that $\rho(D(r)) = g(r)$ for all $r\geq -1$. It only remains to 
prove that $\rho(D(r)) = g(r)$ for all $r\in\Z$.

To achieve this, we take $p = 2,3,4, \ldots$ and consider a subalgebra
$\L_p \subset \L$:
$$\L_p = \Span \left< D(pk) \big| k\in \Z \right> .$$
It is easy to see that the map $\theta_p: \L_p \rightarrow \L$ defined by
$\theta_p (D(pk)) = p D(k)$, is an isomorphism. Thus everything we proved for 
$\L$ is also valid for $\L_p$. This means that there exists a polynomial $g_p(s)$
such that $\rho(D(pr)) = g_p(pr)$ for all $r\geq -1$. Since the values of the
polynomials $g(s)$ and $g_p(s)$ coincide at infinitely many points, we conclude that
$g_p(s) = g(s)$. Taking now $r=-1$ and letting $p=2,3,\ldots $, we get that
$\rho(D(-p)) = g(-p)$. Thus $\rho(D(r)) = g(r)$ for all $r\in\Z$.
Theorem 1 is now proved.

\

\

{\bf 5. Modules for the semidirect product of $\VT$ with a multi-loop algebra.}

\

Let $\dg$ be a finite-dimensional Lie algebra over $\C$. Consider a multi-loop
Lie algebra
$$\wg = \FT \otimes \dg .$$
The Lie algebra $W_n = \VT$ acts in a natural way on the multi-loop algebra, so we can
form the semidirect product
$$ \g = W_n \oplus \wg .$$
We define the category $\J$ consisting of $\g$-modules $J$ satisfying (J1)-(J3) and
also 
$$ \widetilde g \left( f v \right) = f ( \widetilde g v) ,
\; \; \hbox{\rm for \ } 
\widetilde g \in \wg, \; f \in \FT, \; v \in J .\leqno{(J4)}$$

{\bf Theorem 4.} (a) Every $\g$-module in category $\J$ is a polynomial module.

(b) There exists a 1-1 correspondence between indecomposable $\g$-modules in
category $\J$ and pairs $(\lambda, U)$, where $\lambda \in \C^n/ \Z^n$
and $U$ is a finite-dimensional indecomposable module for the semidirect product
$$ \g^+ = W_n^+ \oplus \C[z_1, \ldots, z_n] \otimes \dg .$$
This correspondence is given by the tensor product decomposition $J = \FT \otimes U$,
and $\g$ acts according to (3.5) and
$$\left( \e{s} g \right) \left( \e{m} v \right) 
= \sum\limits_{\beta \geq 0} {s^\beta \over \beta !}
\e{(m+s)} \rho \left( z^\beta g \right) v , \; \; \; g\in\dg, \; v \in U . \eqno{(5.1)} $$

(c) Irreducible $\g$-modules in category $\J$ are in a 1-1 correspondence with
pairs $(\lambda, V)$, where $\lambda \in \C^n/ \Z^n$ and $V$ is a finite-dimensional 
irreducible $\gl (\C) \oplus \dg$-module. The action of $\g$ on $J = \FT \otimes V$
is given by (3.6) and
$$ \left( \e{s} g \right) \left( \e{m} v \right) = \e{(m+s)} (gv), \; \; \; 
g\in\dg, \; v \in V . \eqno{(5.2)} $$

\

% The proof of Theorem 4 is a simple modification of the proofs of Theorems 2 and 3,
% and is left as an exercise to the reader.

Let us outline the proof of Theorem 4. 
In the same way as in the discussion at the beginning of
Section 3, we can show that an indecomposable $\g$-module in category $\J$
is a tensor product $J = \FT \otimes U$ with $\dim U < \infty$, and the action of 
$\g$ given by (3.1) and for $g \in \dg$:
$$\left( \e{s} g \right) \left( \e{m} v \right) = \e{(m+s)} g(s) v ,$$
for some operators $g(s) \in \End U$, $s\in\Z^n$. In order to prove that $J$
is a polynomial module, we have to show that the family of operators $\{ g(s) \}$
has a polynomial dependence on $s$.

From the commutator relation in $\g$,
$$ \left[ {1\over 2\pi i} \e{s} {\d \over \d x^j}, \e{m} g \right] 
= m_j \e{(s+m)} g ,$$
we get that
$$ \left[ D_j (s), g(m) \right] = m_j \left( g(s+m) - g(m) \right) . \eqno{(5.3)}$$
Applying the above equality to $g(1,\ldots,1)$ we see that
$$ g(1+s_1,\ldots,1+s_n) =  \left[ D_j (s), g(1,\ldots,1) \right] + g(1,\ldots,1).$$   
By Theorem 2, the operators $\{ D_j (s) \}$ depend on $s$ polynomially, hence the same
is true for $\{ g(s) \}$.

To prove part (b), we expand the polynomial $g(s)$:
$$ g(s) = \sum\limits_{\beta \geq 0}^{\rm finite}
{s^\beta \over \beta !} g^{(\beta)},$$
with $g^{(\beta)} \in \End U$ and $g^{(\beta)} = 0$ for $\beta \gg 0$.

From the expansions of the relations
$$ [g(s), h(m)] = [g,h] (s+m) , \; \; g,h \in \dg ,$$
and (5.3), we see that 
$$ [g^{(\alpha)}, h^{(\beta)} ] = [g,h]^{(\alpha + \beta)}, \eqno{(5.4)}$$
$$ [D_j^{(\alpha)}, g^{(\beta)} ] = \beta_j g^{(\alpha+\beta-\epsilon_j)} .
\eqno{(5.5)}$$
This shows that every indecomposable $\g$-module $J \in \J$ yields
a finite-dimensional $\g^+$-module $U$.

Using the same technique as in the proof of Theorem 3 (b), we can show that
for any finite-dimensional $\g^+$-module $U$, $\rho \left( z^\alpha g \right) = 0$
for $\alpha \gg 0$. This proves that the correspondence is bijective.

Finally, to prove part (c) of the theorem, we note that irreducible $\g$-modules $J$
correspond to irreducible $\g^+$-modules $V$. Let us show that $V$ is in fact a 
$\gl (\C) \oplus \dg$-module. Let $V_\nu$ be an eigenspace for $\rho(E)$,
$E = z_1 {\d \over \d z_1} + \ldots z_n {\d \over \d z_n}$. It is easy to
see that 
$$\rho\left(z^\alpha {\d \over \d z_j} \right) V_\nu \subset V_{\nu + |\alpha| - 1},
\; \; \hbox{\rm and} \; \;  
\rho\left(z^\beta g \right) V_\nu \subset V_{\nu + |\beta|} .$$
Thus $\mathop\oplus\limits_{k=1}^\infty V_{\nu + k}$ is a $\g^+$-submodule in $V$.
Irreducibility of $V$ implies that $V_{\nu + k} = (0)$ for all $k \geq 1$,
and hence the ideal
$$g^{++} = \Span \left< z^\alpha {\d \over \d z_j} , \; z^\beta g
\bigg| \alpha, \beta \in \Z_+^n, \; |\alpha| > 1, |\beta| \geq 1, \; 
g \in \dg, \; j=1,\ldots,n \right> $$
vanishes in every finite-dimensional irreducible $\g^+$-module $V$. The claim of part (c)
follows from the fact that $\g^+ / \g^{++} \cong \gl (\C) \oplus \dg$.

\

{\bf Acknowledgments.} This research is supported by the  Natural Sciences and
Engineering Research Council of Canada.

I thank Mikhail Kochetov for helpful discussions.
I also thank Eswara Rao for keeping me updated on his research.
In particular, this work was strongly influenced by Eswara Rao's papers 
[R1] and [R2]. 

\

\

{\bf References:}

\

\item{[BB]} S.~Berman, Y.~Billig, 
{\it Irreducible representations for toroidal Lie algebras.}
J.Algebra {\bf 221} (1999), 188-231.

\item{[BZ]} Y.~Billig,  K.~Zhao, {\it Weight modules over exp-polynomial
Lie algebras.} J. Pure Appl. Algebra {\bf 191} (2004), 23-42.

\item{[R1]} S.~Eswara Rao, {\it A new class of modules for derivations of Laurent polynomial ring in $n$ variables.} Preprint, 2003.

\item{[R2]} S.~Eswara Rao, {\it Partial classification of modules 
for Lie algebra of diffeomorphisms of $d$-dimensional torus.}
J. Math. Phys. {\bf 45} (2004), 3322-3333.

\item{[M]} O.~Mathieu, {\it Classification of Harish-Chandra modules 
over the Virasoro algebra.} Invent. Math. {\bf 107} (1992), 225-234. 

\item{[O]} P.J.~Olver, {\it Applications of Lie groups to differential equations.} Graduate Texts in Math. {\bf 107}, Springer-Verlag, New York, 1986.

\item{[Ru]} A.N.~Rudakov, {\it Irreducible representations of
infinite-dimensional Lie algebras of Cartan type.} Math. USSR Izv. 
{\bf 8} (1974), 836-866.

\item{[S]} D.J.~Saunders, {\it The geometry of jet bundles.} London Math. Soc. Lecture Notes Ser. {\bf 142}, Cambridge Univ. Press, Cambridge, 1989.

\item{[vdW]} B.L.~van der Waerden, {\it Modern Algebra.} Frederick Ungar Publ., New York, 1949.

\

\

{\it School of Mathematics and Statistics}

{\it Carleton University}

{\it Ottawa, Ontario,}

{\it K1S 5B6, Canada.}

\end